\documentclass[reqno]{amsart}

\theoremstyle{plain}
\newtheorem{lemma}{Lemma}
\newtheorem{theorem}[lemma]{Theorem}

\theoremstyle{remark}

\def\aa{\alpha}

\def\tt{\theta}

\def\Dd{\Delta}

\def\Om{\Omega}

\def\pp{\partial}

\begin{document}

\title[Regularity of 3D Navier--Stokes Equations]
{Regularity Criterion for Solutions of Three-Dimensional Turbulent
Channel Flows}

\date{September 23, 2006}

\author[C. Cao]{Chongsheng Cao}
\address[C. Cao]
{Department of Mathematics  \\
Florida International University  \\
University Park  \\
Miami, FL 33199, USA}
\email{caoc@fiu.edu}

\author[J. Qin]{Junlin Qin}
\address[J. Qin]
{Department of Mathematics \\
Xian Jiaotong University \\
Xian, Shaanxi, China} \email{jlqin@mail.xjtu.edu.cn}

\author[E.S. Titi]{Edriss S. Titi}
\address[E.S. Titi]
{Department of Mathematics \\
and  Department of Mechanical and  Aerospace Engineering \\
University of California \\
Irvine, CA  92697-3875, USA \\
{\bf also}  \\
Department of Computer Science and Applied Mathematics \\
Weizmann Institute of Science  \\
Rehovot 76100, Israel} \email{etiti@math.uci.edu  and edriss.titi@weizmann.ac.il}

\begin{abstract}
In this paper we consider the three-dimensional  Navier--Stokes
equations in infinite channel. We provide a regularity criterion for
solutions of the three-dimensional Navier--Stokes equations in terms
of the vertical component of the velocity field.
\end{abstract}

\maketitle



\section{Introduction}   \label{S-1}

Turbulence stands out as a prototype of multi-scale phenomenon that
occurs in nature. It involves  wide ranges of spatial and temporal
scales which makes it very difficult to study analytically and
prohibitively expensive to simulate computationally.
Turbulent channel flows  are considered to be the simplest flows
confined within physical boundaries that can be simulated
numerically and that demonstrates many of the common features of
turbulence. In this paper we consider three-dimensional finite
energy turbulent flows of viscous incompressible homogeneous fluids
in the infinite channel $\Om =\mathbb{R}^2 \times [-L, L] \subset
\mathbb{R}^3$, subject to the no-slip Dirichlet boundary conditions.
These flows are governed by the three-dimensional Navier--stokes
system of equations:
\begin{eqnarray}
&&\hskip-0.5in \frac{\partial u}{\partial t} - \nu \Dd u + ( u \cdot
\nabla )u + \nabla p = f     \hskip0.5in \mbox{in  }  \Om
\label{NS1}   \\
&&\hskip-0.5in
 \nabla \cdot u=0    \hskip1.86in  \mbox{in  }  \Om \label{NS2}   \\
&&\hskip-0.5in
 u=0   \hskip2.07in  \mbox{on  }  \partial \Om \label{NS3}   \\
&&\hskip-0.5in
 \lim_{x \to \infty} u(t,x) = 0
\label{boundary-infinity}   \\
&&\hskip-0.5in
 u(0,x) =u_0 (x)    \hskip1.5in  \mbox{in  }  \Om.   \label{NS4}
\end{eqnarray}
Here, $u=(u_1, u_2, u_3)$ represents the unknown velocity vector
field, and $p$ is the unknown pressure scalar; where $\nu >0$,  the
constant kinematic viscosity,  $f$,  the body forcing term, and
$u_0$, the initial velocity, are given.

Mathematically, it is well-known that the three-dimensional system
(\ref{NS1})--(\ref{NS4}) has  global (for all time and all initial
data) weak solutions (see, e.g., \cite{CF88}, \cite{GA94},
\cite{LADY1969}, \cite{Sohr}, \cite{TT83}, \cite{TT84} and
references therein). The question of well-posedness, in the sense of
Hadamard, and in particular the question of uniqueness, of these
weak solutions is still an open problem. On the other hand, it is
also well-established (see, e.g., \cite{CF88}, \cite{GA94},
\cite{LADY1969}, \cite{Sohr}, \cite{TT83}, \cite{TT84} and
references therein) that the system (\ref{NS1})--(\ref{NS4})
possesses a unique strong (regular) solution, which depends
continuously on the initial data, for a short interval of time
$[0,T_*)$, where $T_*$ depends on the size of initial datum, $u_0$,
on $f,\nu$ and $L$. Moreover, it is also well-known that the
existence (for all time) and uniqueness of strong (regular)
solutions is guaranteed under suitable additional assumptions (see,
e.g., \cite{BL02}, \cite{BG02},
\cite{CL01},\cite{Cao-Titi-Pressure}, \cite{HE02}, \cite{KT84},
\cite{KZ05}, \cite{PO03}, \cite{PG59}, \cite{SJ62}, \cite{ZY02} and
references therein). In particular, some of these recent results
involve conditions on only one component of the velocity field of
the $3D$ NSE in the whole space $\mathbb{R}^3$ or under periodic
boundary conditions (see, e.g., \cite{HE02}, \cite{KZ05},
\cite{PO03}, \cite{ZY02}). In this paper, we study this type of
sufficient conditions for the global regularity of the $3D$ NSE in
the infinite channel $\Omega$, subject to no-slip Dirichlet boundary
condition on the physical boundary of the channel. Using the
geophysical terminology, our condition is formulated in terms of the
third component of the baroclinic mode $\widetilde{u}_3$ (see
(\ref{AVG}), below, for the definition of the barotropic mode
(vertically averaged mode), $\bar{u}$, and the baroclinic mode (the
fluctuation about the barotropic mode)). Specifically, our results
states that if $\widetilde{u}_3$ satisfies
\begin{equation}
\nabla \widetilde{u}_3 \in L^{\infty} ([0, \infty), L^2(\Om)),
\label{CON}
\end{equation}
then the strong (regular) unique solution of the 3D Naiver-Stokes
equation (\ref{NS1})--(\ref{NS4}) exists for all time.

Let us observe that our condition (\ref{CON}) seems to be slightly
tighter than the former ones (cf. e.g., \cite{HE02}, \cite{KZ05},
\cite{PO03}, \cite{ZY02}). However, unlike the previous works we
study here the $3D$ Navier--Stokes in a domain with physical
boundaries under the no-slip  Dirichlet boundary conditions.
Furthermore, we emphasize that the techniques developed here, which
are inspired by ideas presented in \cite{CT07}, are totally
different than the previous ones.

\vskip0.1in

Let us denote by $L^q(\Om), L^q(\mathbb{R}^2)$, and $H^m (\Om), H^m (\mathbb{R}^2)$ the
usual $L^q-$Lebesgue and Sobolev spaces, respectively
(\cite{AR75}). We denote by
\begin{equation}
\| \phi\|_q = \left\{ \begin{array}{ll} \left(  \int_{\Om}
|\phi|^q \; dx_1dx_2dx_3
\right)^{\frac{1}{q}},   \qquad & \mbox{ for every $\phi \in L^q(\Om)$} \\
\left(  \int_{\mathbb{R}^2}  |\phi |^q \; dx_1dx_2 \right)^{\frac{1}{q}}, &
\mbox{ for every $\phi \in L^q(\mathbb{R}^2)$}.
\end{array} \right.
 \label{LQ}
\end{equation}
Let
\begin{eqnarray*}
\mathcal{V} &=&  \left\{ v \in C^{\infty}_0(\Om): \nabla \cdot v =
0 \right\}.
\end{eqnarray*}
Since we are interested in flows of finite energy in the infinite
channel $\Omega$, we consider the spaces $H$ and $V$, defined to be
the closures of the set $\mathcal{V}$ in $L^2(\Om)$ under
$L^2-$topology, and  in $H^1(\Om)$ under $H^1-$topology,
respectively.
 Denote by $P:L^2 \rightarrow H$,
the orthogonal projection, and let   $A=-P\Dd$ be the Stokes
operator subject to the homogeneous Dirichlet  boundary condition
(\ref{NS3}). It is well known that the Navier--Stokes equations
(NSE) (\ref{NS1})--(\ref{NS4}) are equivalent to the functional
differential equation (see, e.g.,  \cite{CF88}, \cite{Sohr},
\cite{TT83},\cite{TT84})
\begin{eqnarray}
&&\hskip-0.5in
\frac{du}{dt} + \nu A u +B(u,u) = f, \label{EQ}   \\
&&\hskip-0.5in
 u(0) =u_0,    \label{EQ1}
\end{eqnarray}
where $B(u,u) = P ( (u \cdot \nabla ) u ),$ the nonlinear (bilinear)
term. We say $u$ is a Leray--Hopf weak solution to the system
(\ref{EQ})--(\ref{EQ1}) if $u$ satisfies (see, e.g.,  \cite{CF88},
\cite{TT83},\cite{TT84})

\begin{itemize}

\item[(1)]
$u  \in C([0, T], H\text{-weak}) \cap L^2([0, T], V),  $ and $ \pp_t
u \in L^1([0, T], V^{\prime}),$ where $V^{\prime}$ is the dual space
of $V$,

\item[(2)] the weak formulation:
\begin{eqnarray*}
&&\hskip-0.35in \int_{\Om} u (t,x)\cdot \phi (x) \, dx
-\int_{\Om}  u(t_0,x) \cdot \phi(x)  \, dx  \nonumber  \\
&&\hskip-0.25in =  - \int_{t_0}^t \int_{\Om}  \left(\nu (\nabla
u(s,x) : \nabla  \phi(x)) + ( u(s,x) \cdot \nabla)  u(s,x) \cdot
\phi(x) \right) \, dx \, ds +\int_{t_0}^t \int_{\Om} \left (  f(s,x)
\cdot \phi (x) \right) \, dx\, ds, \label{WEAK}
\end{eqnarray*}
for every  $\phi \in \mathcal{V},$ and
almost every $t$, $t_0\in [0,T]$.

\item[(3)] the energy inequality:
\begin{eqnarray}
&&\hskip-.68in \|u(t)\|_2^2 -  \|u(t_0)\|_2^2 +2 \nu \int_{t_0}^t \|
\nabla u (s) \|_2^2 \; ds
 \leq  2 \int_{t_0}^t\int_{\Om} f(s,x) \cdot u(s,x) \; dxds,  \label{ENG}
\end{eqnarray}
\end{itemize}
for all $t\in[0,T]$, and for almost every $t_0$ in the interval
$[0,t]$.

Moreover, a weak solution is called strong solution of
(\ref{EQ})--(\ref{EQ1}) on $[0,T]$ if, in addition, it satisfies
\begin{eqnarray*}
&&u  \in C([0,T], V) \cap L^2([0,T], H^2 (\Om)).
\end{eqnarray*}

 \vskip0.1in

For convenience, we recall the following Gagiliardo-Nirenberg,
Ladyzhenskaya, and Sobolev inequalities  (cf. e.g., \cite{AR75},
\cite{CF88}, \cite{GA94} \cite{LADY}, \cite{LADY1969} and
\cite{Sohr}) in $\mathbb{R}^2$:
\begin{eqnarray}
&&\hskip-.68in \| \phi\|_{L^r (\mathbb{R}^2)} \leq C_r \| \phi
\|_{L^2(\mathbb{R}^2)}^{2/r} \| \phi
\|_{H^1(\mathbb{R}^2)}^{\frac{r-2}{r}}, \qquad  r < \infty, \label{SI-2}
\end{eqnarray}
for every $\phi \in H^1(\mathbb{R}^2),$ and in $\mathbb{R}^3$:
\begin{eqnarray}
&&\hskip-.68in \| \psi \|_{L^{\aa}(\Om)} \leq C_{\aa} \| \psi
\|_{L^2(\Om)}^{\frac{6-\aa}{2\aa}} \| \psi
\|_{H^1(\Om)}^{\frac{3({\aa}-2)}{2{\aa}}},  \label{SI1}
\end{eqnarray}
for every $u\in H^1(\Om), 2 \leq \aa \leq 6.$  Here $C_r$ and
$C_{\aa}$ are scale invariant  constants. We also recall the
Pioncar\'{e} inequality:
\begin{eqnarray}
&&\hskip-.7in \|\nabla v\|_2 \geq \frac{C_0}{L} \|v\|_2   \hskip2.0in  \forall  v
\in V
\label{P-1}  \\
&& \hskip-.7in \|Av\|_2  \geq \frac{C_0}{L}  \|\nabla v\|_2 \hskip2.0in  \forall v
\in {\mathcal{D}}(A), \label{P-2}
\end{eqnarray}
where $C_0$ is a scale invariant constant.
 Also, we recall the integral version of
Minkowsky inequality for the $L^r$ spaces, $r\geq 1$. Let $\Om_1
\subset \mathbb{R}^{m_1}$ and
 $\Om_2 \subset \mathbb{R}^{m_2}$ be two measurable sets, where
$m_1$ and $m_2$ are two positive integers. Suppose that
$\phi(\xi,\eta)$ is measurable over $\Om_1 \times \Om_2$. Then,
\begin{equation}
\hskip0.35in \left[ { \int_{\Om_1} \left( \int_{\Om_2}
|\phi(\xi,\eta)| d\eta \right)^r d\xi } \right]^{1/r} \leq
\int_{\Om_2} \left( \int_{\Om_1} |\phi(\xi,\eta)|^r d\xi
\right)^{1/r} d\eta. \label{MKY}
\end{equation}


\section{Global Existence of the Strong Solution} \label{S-2}

In this section we will show the global existence of the strong
solutions to the three-dimensional Navier--Stokes system
(\ref{NS1})--(\ref{NS4}) under assumption (\ref{CON}).

We will denote by
\begin{equation}
\hskip0.35in \overline{\tt} (x_1, x_2) = \frac{1}{2L} \; \int_{-L}^L \tt (x_1, x_2, x_3) \; dx_3
\qquad  \mbox{and}  \quad \widetilde{\tt} = \tt - \overline{\tt}. \label{AVG}
\end{equation}
Following the geophysical fluid dynamics terminology we will call
$\overline{\tt}$   the barotropic mode and $\widetilde{\tt}$ the
baroclinic mode.

{F}rom now on, we will denote by $\nabla_h = (\frac{\pp}{\pp x_1},
\frac{\pp}{\pp x_2})$ and $\Dd_h = \frac{\pp^2}{\pp x_1^2}+
\frac{\pp^2}{\pp x_2^2}$. First, let us prove the following Lemma.

\begin{lemma} \label{LLL}
Suppose that $\xi (x_1, x_2) \in H^1(\mathbb{R}^2),
\phi \in H^1(\Om)$ and $\psi \in L^2(\Om)$.
Then,
\begin{eqnarray*}
&&\hskip-.65in
 \int_{\Om} |\xi| \, |\phi| \, |\psi| \;
 dx_1dx_2dx_3   \\
&&\hskip-.6in
\leq C \|\xi \|_2^{1/2} \; \left( \|\xi \|_2
+ \| \nabla_h \xi \|_2  \right)^{1/2}
 \, \|\phi\|_2^{1/2} \, \left( \|\phi \|_2
+ \| \nabla_h \phi \|_2  \right)^{1/2} \; \|\psi\|_2.
\end{eqnarray*}
\end{lemma}

\begin{proof}
Notice that
\begin{eqnarray*}
&&\hskip-.65in
 \int_{\Om} |\xi| \, |\phi| \, |\psi| \;
 dx_1dx_2dx_3 = \int_{\mathbb{R}^2} \left[ |\xi| \,\left( \int_{-L}^L |\phi| \, |\psi|
 \;dx_3 \right) \right]\; dx_1dx_2.
\end{eqnarray*}
We  will estimate the above term by applying the same method used to establish
Proposition 2.2 in \cite{CT03}. First, by Cauchy--Schwarz
inequality, we obtain
\begin{eqnarray*}
&& \int_{-L}^L |\phi| \, |\psi| \; dx_3 \leq  \left( \int_{-L}^L
|\phi|^2 \; dx_3 \right)^{\frac{1}{2}} \; \left( \int_{-L}^L |\psi|^2
\; dx_3\right)^{\frac{1}{2}}.
\end{eqnarray*}
Thus, by the above and H\"{o}lder inequality, we reach
\begin{eqnarray*}
&& \int_{\Om} |\xi| \, |\phi| \, |\psi| \;
 dx_1dx_2dx_3 \leq \int_{\mathbb{R}^2} \left[ |\xi|\;\left( \int_{-L}^L |\phi|^2
\; dx_3 \right)^{\frac{1}{2}} \; \left( \int_{-L}^L |\psi|^2 \;
dx_3\right)^{\frac{1}{2}}  \; \right] dx_1dx_2  \\
&& \leq  \left(  \int_{\mathbb{R}^2} |\xi|^4 dx_1dx_2
\right)^{\frac{1}{4}} \left[ { \int_{\mathbb{R}^2} \left( {
\int_{-L}^L \left| \phi \right|^2 dx_3 }\right)^2 dx_1dx_2
}\right]^{\frac{1}{4}} \left[ { \int_{\Om} \left| \psi  \right|^2
dx } \right]^{\frac{1}{2}}.
\end{eqnarray*}
By using Minkowsky  inequality (\ref{MKY}),  we get
\[
\left[ { \int_{\mathbb{R}^2} \left( { \int_{-L}^L \left| \phi
\right|^2 dx_3 }\right)^2 dx_1dx_2 }\right]^{\frac{1}{2}} \leq
\int_{-L}^L \left( { \int_{\mathbb{R}^2} \left| \phi \right|^4
dx_1dx_2 }\right)^{\frac{1}{2}} dx_3.
\]
Thanks to (\ref{SI-2}) with $r=4$, for every fixed $x_3$ we have
\begin{eqnarray*}
&& \left({ \int_{\mathbb{R}^2}   \left| \phi \right|^4 dx_1dx_2
}\right)^{\frac{1}{4}}   \leq C \left\| \phi
\right\|_{L^2({\mathbb{R}^2})}^{\frac{1}{2}} \left\| \phi
\right\|_{H^1({\mathbb{R}^2})}^{\frac{1}{2}}.
\end{eqnarray*}
As a result of the above and the Cauchy--Schwarz inequality, we
obtain
\begin{eqnarray*}
&& \int_{-L}^L \left( {  \int_{\mathbb{R}^2}   \left| \phi \right|^4
dx_1dx_2
}\right)^{\frac{1}{2}} dx_3 \\
&& \leq C \int_{-L}^L  \left\| \phi \right\|_{L^2({\mathbb{R}^2})}
\left\| \phi \right\|_{H^1({\mathbb{R}^2})} dx_3  \\
&&
 \leq C \left( {
 \int_{-L}^L \left\| \phi \right\|_{L^2({\mathbb{R}^2})}^2 \; dx_3
} \right)^{\frac{1}{2}} \; \left( {
 \int_{-L}^L
\left\| \phi \right\|_{H^1({\mathbb{R}^2})}^2 \; dx_3
} \right)^{\frac{1}{2}}  \\
&& \leq  C \left( \| \phi  \|_2 +  \| \nabla_h \phi \|_2 \right) \; \| \phi  \|_2.
\end{eqnarray*}
Therefore,
\begin{eqnarray}
&&  \left[ { \int_{\mathbb{R}^2} \left( { \int_{-L}^L  \left| \phi
\right|^2 dx_3 }\right)^2 dx_1dx_2
}\right]^{\frac{1}{4}}
\nonumber  \\
&& \leq C \left( \| \phi  \|_2 +  \| \nabla_h \phi \|_2 \right)^{1/2} \; \| \phi  \|_2^{1/2}.
\label{NL-1}
\end{eqnarray}
By using (\ref{SI-2}) with $r=4$, we have
\begin{eqnarray*}
&& \left[ { \int_{\mathbb{R}^2}  \left| \xi \right|^4
dx_1dx_2 } \right]^{\frac{1}{4}}   \leq C
\|\xi\|_{L^2({\mathbb{R}^2})}^{1/2}
\| \xi\|_{H^1({\mathbb{R}^2})}^{1/2}.
\end{eqnarray*}
Thus, by the above and (\ref{NL-1}), we get
\begin{eqnarray*}
&&\hskip-.65in
 \int_{\Om} |\xi| \, |\phi| \, |\psi| \;
 dx_1dx_2dx_3   \\
&&\hskip-.6in
\leq C \|\xi \|_2^{1/2} \; \left( \|\xi \|_2
+ \| \nabla_h \xi \|_2  \right)^{1/2}
 \, \|\phi\|_2^{1/2} \, \left( \|\phi \|_2
+ \| \nabla_h \phi \|_2  \right)^{1/2} \; \|\psi\|_2.
\end{eqnarray*}

\end{proof}

\begin{theorem} \label{T-MAIN}
Let $f \in L^{\infty}([0, \infty), L^{2}(\Om)), u_0 \in V$. Let
$u=(u_1,u_2,u_3)$ be a weak solution of the system {\em
(\ref{EQ})--(\ref{EQ1})} in $[0,\infty)$. Suppose that for $T>0,$
$\nabla \widetilde{u}_3 \in L^{\infty}([0, T], L^{2}(\Om))$; that
is, $\widetilde{u}_3$ satisfies
\begin{eqnarray}
&&\hskip-.65in
 \displaystyle{\sup_{0\leq t \leq T} } \| \nabla \widetilde{u}_3 (t)\|_2  < \infty.
\label{SUP}
\end{eqnarray}
Then $u$ is the strong solution of the system {\em
(\ref{EQ})--(\ref{EQ1})} on $[0, T]$.
\end{theorem}

\vskip0.05in

\begin{proof} Let $u_0\in V$. Following, for instance, the
Galerkin method one can show that there exists a unique strong
solution $u$ for the system ({\ref{EQ}})--(\ref{EQ1}), with the
initial datum $u_0$, for a short interval of time (see, e.g.,
\cite{CF88}, \cite{LADY1969}, \cite{Sohr},\cite{TT83} and
\cite{TT84}. Suppose that $[0,T_*)$ is the maximal interval of
existence of this strong solution $u$. It is also well known (see,
e.g., the above references) that there exists a Leray-Hopf weak
solution for the system ({\ref{EQ}})--(\ref{EQ1}), with the same
initial datum $u_0$, which exists globally in time, i.e. for all
time $t\ge 0$. Most importantly, following the work of J.~Sather and
J.~Serrin in \cite{Serrin-Uniqueness} on can show that all the
Leray-Hopf weak solutions coincide with the unique strong solution,
$u$, on the interval $[0,T_*)$ .

To conclude our proof we need to show that $T<T_*$. Suppose, arguing
by contradiction, that $T_* \le T$, and that (\ref{SUP}) holds. If
we show that the $\limsup_{t \to T_*^-} \|u(t)\|_{H^1} < \infty$
then $[0,T_*)$ is not a maximal interval of existence, which leads
to a contradiction.

For the rest of this proof we consider the strong solution, $u$, in
the the interval $[0,T_*)$. From the energy inequality of
 (\ref{ENG}), which is satisfied by all the Leray-Hopf weak solutions (see, for example,
\cite{CF88}, \cite{LADY1969}, \cite{Sohr}, \cite{TT83} or
\cite{TT84} for details), we have
\begin{eqnarray}
&&\hskip-.68in \|u(t)\|_2^2 \leq C \frac{L^4\, F^2}{\nu^2} +
e^{-\frac{\nu \, t}{L^2}}
\|u_0\|_2^2, \label{L2}  \\
&&\hskip-.68in  \nu \int_0^t \| \nabla u (s) \|_2^2 \; ds \leq C
\frac{L^2\, F^2 \,  t}{\nu} + \|u_0\|_2^2, \label{L22}
\end{eqnarray}
for all  $t\in [0,T_*)$, where
\begin{equation}
F =\|f\|_{L^\infty([0,\infty),L^2(\Om))}.\label{F}
\end{equation}
In particular, since $t < T_* \leq T$ we have
\begin{eqnarray}
&&\hskip-.68in \|u(t)\|_2^2 + \nu \int_0^t \| \nabla u (s) \|_2^2
\; ds \leq K_1, \label{K1}
\end{eqnarray}
where
\begin{eqnarray}
&&\hskip-.68in K_1= C \frac{F^2\, (L^4+\nu \, T)}{\nu^2} + 2
\|u_0\|_2^2. \label{K-1}
\end{eqnarray}
Taking the  inner product of the  equation (\ref{EQ}) with
$-\Dd_h u$ in $H$, and notice that $P \pp_{x_i} = \pp_{x_i} P$ for $i=1,2$,  we get
\begin{eqnarray*}
&&\hskip-.68in \frac{1}{2} \frac{d \|\nabla_h u \|_2^2 }{d t} + \nu
\left\|\Dd_h u \right\|_2^2 + \nu
\left\|\nabla_h u_z \right\|_2^2  =  - \int_{\Om} \left(f - B(u, u)\right) \cdot
\Dd_h u \;  dx_1dx_2dx_3.
\end{eqnarray*}
By integration by parts we get
\begin{eqnarray*}
&&\hskip-.265in    -\int_{\Om} B(u, u) \cdot \Dd_h u \; dx_1dx_2dx_3 =
\int_{\Om} \sum_{l=1}^2 \frac{\pp u_k}{\pp x_j}\frac{\pp u_j}{\pp x_l}\frac{\pp
u_k}{\pp x_l} \; dx_1dx_2dx_3   \\
&&\hskip-.265in =  \int_{\Om} \left\{ \left(\frac{\pp u_1}{\pp
x_1}\right)^3+\left(\frac{\pp u_2}{\pp x_2}\right)^3 +  \left(\frac{\pp u_1}{\pp x_1}+\frac{\pp u_2}{\pp
x_2}\right)  \left[ \left(\frac{\pp u_1}{\pp x_2}\right)^2 +\left(\frac{\pp u_2}{\pp x_1}\right)^2
+ \frac{\pp u_1}{\pp x_2}\;\frac{\pp u_2}{\pp x_1} \right]  \right.  \\
&&\hskip-.0015in
 +  \left. \sum_{k, l=1}^2 \frac{\pp u_k}{\pp x_3}\frac{\pp
u_3}{\pp x_l}\frac{\pp u_k}{\pp x_l} + \sum_{j,l=1}^2  \frac{\pp u_3}{\pp x_j}\frac{\pp u_j}{\pp x_l}\frac{\pp
u_3}{\pp x_l} + \sum_{l=1}^2  \frac{\pp u_3}{\pp x_3}\frac{\pp u_3}{\pp x_l}\frac{\pp
u_3}{\pp x_l}  \right\} \; dx_1dx_2dx_3    \\
&&\hskip-.265in =    \int_{\Om} \left\{ - \; \frac{\pp u_3}{\pp x_3} \left[
\left(\frac{\pp u_1}{\pp
x_1}\right)^2+\left(\frac{\pp u_2}{\pp x_2}\right)^2 -\frac{\pp u_1}{\pp x_1}\frac{\pp u_2}{\pp x_2}
+  \left(\frac{\pp u_1}{\pp x_2}\right)^2 +\left(\frac{\pp u_2}{\pp x_1}\right)^2 +
\frac{\pp u_1}{\pp x_2}\;\frac{\pp u_2}{\pp x_1}
 \right] \right.   \\
&&\hskip-.0015in +   \sum_{k, l=1}^2 \left(  \frac{\pp u_k}{\pp x_l} \frac{\pp
u_3}{\pp x_3}\frac{\pp u_k}{\pp x_l}  + u_k \frac{\pp
u_3}{\pp x_3} \frac{\pp^2
u_k}{\pp x_l \pp x_l}  - u_k  \frac{\pp
u_3}{\pp x_l}\frac{\pp^2 u_k}{\pp x_l  \pp x_3}   \right)    \\
&&\hskip-.0015in +   \left.
\sum_{j,l=1}^2  \frac{\pp u_3}{\pp x_j}\frac{\pp u_j}{\pp x_l}
\frac{\pp u_3}{\pp x_l} - \sum_{l=1}^2  \left(\frac{\pp u_1}{\pp x_1}+\frac{\pp u_2}{\pp
x_2}\right)  \frac{\pp u_3}{\pp x_l}\frac{\pp
u_3}{\pp x_l}  \right\} \; dx_1dx_2dx_3
 \\
&&\hskip-.265in =    \int_{\Om} \left\{ - \frac{\pp
\widetilde{u}_3}{\pp x_3} \left[ \left(\frac{\pp u_1}{\pp
x_1}\right)^2+\left(\frac{\pp u_2}{\pp x_2}\right)^2 -\frac{\pp
u_1}{\pp x_1}\frac{\pp u_2}{\pp x_2} +  \left(\frac{\pp u_1}{\pp
x_2}\right)^2 +\left(\frac{\pp u_2}{\pp x_1}\right)^2 + \frac{\pp
u_1}{\pp x_2}\;\frac{\pp u_2}{\pp x_1}
 \right] \right.  \\
&&\hskip-.0015in  +  \sum_{k, l=1}^2 \left(  \frac{\pp u_k}{\pp x_l}
\frac{\pp \widetilde{u}_3}{\pp x_3}\frac{\pp u_k}{\pp x_l}  + u_k
\frac{\pp \widetilde{u}_3}{\pp x_3} \frac{\pp^2 u_k}{\pp x_l \pp
x_3} - u_k  \frac{\pp \widetilde{u}_3}{\pp x_l}\frac{\pp^2 u_k}{\pp
x_l \pp x_3}
\right)    \\
&&\hskip-.0015in - \overline{u}_3 \left[ \sum_{j,l=1}^2
 \left[ \frac{\pp}{\pp x_j} \left( \frac{\pp u_j}{\pp x_l}
\frac{\pp u_3}{\pp x_l}  \right)
+ \frac{\pp}{\pp x_l} \left( \frac{\pp u_j }{\pp x_3} \frac{\pp u_j}{\pp x_l} \right) \right]
  - \sum_{l=1}^2 \frac{\pp }{\pp x_l} \left[ \left(\frac{\pp u_1}{\pp x_1}+\frac{\pp u_2}{\pp
x_2}\right)  \frac{\pp
u_3}{\pp x_l} \right] \right]       \\
&&\hskip-.0015in   \left. + \sum_{j,l=1}^2  \frac{\pp
\widetilde{u}_3}{\pp x_j}\frac{\pp u_j}{\pp x_l} \frac{\pp u_3}{\pp
x_l} - \sum_{l=1}^2  \left(\frac{\pp u_1}{\pp x_1}+\frac{\pp
u_2}{\pp x_2}\right)  \frac{\pp \widetilde{u}_3}{\pp x_l}\frac{\pp
u_3}{\pp x_l}  \;\right\}  dx_1dx_2dx_3.
\end{eqnarray*}
Then, from the above and the Cauchy--Schwarz inequality we obtain
\begin{eqnarray*}
&&\hskip-.268in \frac{1}{2} \frac{d \|\nabla_h u \|^2 }{d t} + \nu
\|\Dd_h u\|_2^2 + \nu \|\nabla_h u_z\|_2^2  \leq F \;\|\Dd_h u\|_2
+ C \int_{\Om} |\overline{u}_3| \; |\nabla u| \; |\nabla_h \nabla u| \; dx_1dx_2dx_3  \\
&&\hskip-.15in + C \int_{\Om} |\nabla \widetilde{u}_3| \; |\nabla_h
u|^2  \; dx_1dx_2dx_3 + C \int_{\Om} | u| \; |\nabla
\widetilde{u}_3| \; |\nabla_h \nabla u| \; dx_1dx_2dx_3.
\end{eqnarray*}
where $F$ is given in (\ref{F}). By applying Lemma \ref{LLL}  we
obtain
\begin{eqnarray}
&&\hskip-.68in
\int_{\Om} |\overline{u}_3| \; |\nabla u| \; |\nabla_h \nabla u| \; dx_1dx_2dx_3   \nonumber   \\
&&\hskip-.68in \leq C \|\overline{u}_3\|_2^{1/2} \; \left(
\|\overline{u}_3\|_2 + \|\nabla_h \overline{u}_3\|_2 \right)^{1/2}
\;
 \|\nabla u\|_2^{1/2}  \left( \|\nabla u\|_2+  \|\nabla_h \nabla u \|_2 \right)^{1/2}
\|\nabla_h \nabla u\|_2
\nonumber   \\
&&\hskip-.68in \leq C \|\overline{u}_3\|_2 \;\|\nabla u\|_2
\;\|\nabla_h \nabla u\|_2 + C \|\overline{u}_3\|_2^{1/2} \;
\|\nabla_h \overline{u}_3\|_2^{1/2} \;
 \|\nabla u\|_2^{1/2}  \|\nabla_h \nabla u \|_2^{1/2}
\|\nabla_h \nabla u\|_2     \nonumber   \\
&&\hskip-.68in
\leq C \|u\|_2 \;\|\nabla u\|_2 \;\|\nabla_h \nabla u\|_2 +
C \|u\|_2^{1/2} \; \|\nabla_h u \|_2^{1/2} \;
 \|\nabla u\|_2^{1/2}  \|\nabla_h \nabla u \|_2^{3/2}.   \label{EEE_1}
\end{eqnarray}
By H\"{o}lder inequality, we reach
\begin{eqnarray}
&&\hskip-.268in
\int_{\Om} |\nabla \widetilde{u}_3| \; |\nabla_h u|^2  \; dx_1dx_2dx_3    \nonumber   \\
&&\hskip-.268in \leq C \|\nabla \widetilde{u}_3\|_2  \; \|\nabla_h
u\|_4^2 \leq C \|\nabla \widetilde{u}_3\|_2  \; \|\nabla_h
u\|_2^{1/2} \left( \|\nabla_h u\|_2+  \|\nabla_h \nabla u \|_2
\right)^{3/2}. \label{EEE_2}
\end{eqnarray}
And also by H\"{o}lder inequality, we obtain
\begin{eqnarray}
&&\hskip-.268in
\int_{\Om} | u| \; |\nabla \widetilde{u}_3| \; |\nabla_h \nabla u| \; dx_1dx_2dx_3  \nonumber   \\
&&\hskip-.268in \leq C \int_{\mathbb{R}^2} \left\{  \|u(x_1, x_2,
\cdot)\|_{\infty} \|\nabla \widetilde{u}_3 (x_1, x_2, \cdot) \|_2
\; \|\nabla_h \nabla u (x_1, x_2, \cdot)\|_2 \right\} \; dx_1 dx_2  \nonumber   \\
&&\hskip-.268in \leq C \int_{\mathbb{R}^2} \left\{  \|u(x_1, x_2,
\cdot)\|_2^{1/2} \|\frac{\pp u}{\pp x_3} (x_1, x_2, \cdot)\|_2^{1/2}
\|\nabla \widetilde{u}_3 (x_1, x_2, \cdot) \|_2
\; \|\nabla_h \nabla u (x_1, x_2, \cdot)\|_2 \right\} \; dx_1 dx_2  \nonumber   \\
&&\hskip-.268in
\leq C \left\{ \int_{\mathbb{R}^2}  \|u (x_1, x_2, \cdot)\|_2^{4}  \; dx_1 dx_2  \right\}^{1/8}
\left\{ \int_{\mathbb{R}^2}  \|\frac{\pp u}{\pp x_3} (x_1, x_2, \cdot) \|_2^{4}  \; dx_1 dx_2  \right\}^{1/8}
\nonumber   \\
&&\hskip-.158in \times \left\{ \int_{\mathbb{R}^2}  \|\nabla
\widetilde{u}_3 (x_1, x_2, \cdot) \|_2^{4}  \; dx_1 dx_2
\right\}^{1/4}
\|\nabla_h \nabla u \|_2   \nonumber   \\
&&\hskip-.268in \leq C\|u \|_2^{1/4}  \|\nabla_h u \|_2^{1/4} \;
\left\|\frac{\pp u}{\pp x_3}  \right\|_2^{1/4} \; \left\|\frac{\pp
\nabla_h u}{\pp x_3}  \right\|_2^{1/4} \;  \left\|\nabla
\widetilde{u}_3  \right\|_2^{1/2} \; \|\nabla_h \nabla
\widetilde{u}_3 \|_2^{1/2} \; \|\nabla_h \nabla u \|_2.
\label{EEE_3}
\end{eqnarray}
By (\ref{EEE_1})--(\ref{EEE_3}) and Young's inequality we get
\begin{eqnarray*}
&&\hskip-.68in  \frac{d \|\nabla_h u \|_2^2 }{d t} +
\nu \|\nabla_h \nabla u\|_2^2
\leq C F^2+ C \|u\|_2^{2}  \|\nabla u\|_2^{2}   \\
&&\hskip-.158in
+ C \left(  \|u\|_2^{1/2} \;  \|\nabla u\|_2^{2}  +
  \|\nabla \widetilde{u}_3\|_2^4 +
\|u \|_2^{2}   \; \left\|\frac{\pp u}{\pp x_3}  \right\|_2^{2} \;
 \left\|\nabla \widetilde{u}_3  \right\|_2^{4}
\right) \|\nabla_h u \|_2^{2},
\end{eqnarray*}
where $F$ is given in (\ref{F}). Thanks to Gronwall inequality, we
obtain, for all $t\in[0,T_*)$,
\begin{eqnarray*}
&&\hskip-.68in \|\nabla_h u (t)\|_2^2  + \nu \int_0^t  \|\nabla_h
\nabla u(s)\|_2^2 \; ds \leq K_2,    \label{K-2}
\end{eqnarray*}
where
\begin{eqnarray}
&&\hskip-.68in K_2 = e^{  CK_1^2 + C \left(  T  +K_1^2  \right)
\displaystyle{ \max_{0\leq s \leq T} } \|\nabla \widetilde{u}_3 (s)
\|_2^4 }  \left[ \|u_0\|_{H^1(\Om)}^2  + F^2 + C  K_1^2 \right].
\label{K2}
\end{eqnarray}
 Recall that $\| u \|_V^2 = \int u \cdot Au \; dx_1dx_2dx_3$. It is
well know that $\| u \|_V^2$ is equivalent to $\|\nabla u\|_2^2$
(see, e.g., {\bf \cite{CF88}}). Taking the inner product of the
equation (\ref{EQ}) with $Au$ in $H$,  we get
\begin{eqnarray*}
&&\hskip-.68in \frac{1}{2} \frac{d \| u \|_V^2 }{d t} + \nu
\left\|A u \right\|_2^2  =   \int_{\Om} \left(f - B(u, u)\right) \cdot
A u \;  dx_1dx_2dx_3   \\
&&\hskip-.68in \leq F \|Au\|_2 + C \|u\|_6 \; \|\nabla u \|_3 \; \|Au\|  \\
&&\hskip-.68in \leq F \|Au\|_2 + C \|\nabla_h u\|_2^{2/3}\, \|
\frac{\pp u}{\pp x_3}\|_2^{1/3} \; \|\nabla u \|_2^{1/2} \;
\|Au\|^{3/2}.
\end{eqnarray*}
Here, we used (cf., e.g., \cite{GA94}, p. 33)
\begin{eqnarray*}
&&\hskip-.265in    \|u\|_6 \leq C \left\| \frac{\pp u}{\pp x_1} \right\|_2^{1/3} \;
\left\| \frac{\pp u}{\pp x_2} \right\|_2^{1/3} \; \left\| \frac{\pp u}{\pp x_3} \right\|_2^{1/3}.
\end{eqnarray*}
By Young's inequality we obtain
\begin{eqnarray*}
&&\hskip-.68in
 \frac{d \| u \|_V^2 }{d t} + \nu
\left\|A u \right\|_2^2
 \leq C F^2  + C \|\nabla_h u\|_2^{8/3}\, \| \frac{\pp u}{\pp x_3}\|_2^{4/3} \;
\|u \|_V^2.
\end{eqnarray*}
By Gronwall inequality and (\ref{K-2}), we obtain, for all
$t\in[0,T_*)$,
\begin{eqnarray*}
&&\hskip-.68in \|u (t)\|_V^2  + \nu \int_0^t  \|A
u(s)\|_2^2 \; ds \leq K,    \label{K-F}
\end{eqnarray*}
where
\begin{eqnarray}
&&\hskip-.68in K = e^{  CK_2^{4/3} K_1^{2/3} }  \left[
\|u_0\|_{H^1(\Om)}^2  + F^2  \right].      \label{KF}
\end{eqnarray}
Therefore,
\[
\limsup_{t \to T_*^-} \|u (t)\|_V^2 \le K,
\]
which leads to a contradiction that $[0,T_*)$ is the maximal
interval of existence, and this completes the proof.

\end{proof}

\noindent
\section*{Acknowledgements}
E.S.T.~would like to thank the Bernoulli Center of the  \'{E}cole
Polytechnique F\'{e}d\'{e}ral de Lausanne for the kind hospitality
where this work was completed.  This work was supported in part by
the NSF grant no.~DMS-0504619, the BSF grant no.~2004271, the ISF
grant no.~120/06, and by the MAOF Fellowship of the Israeli Council
of Higher Education.

\end{document}